\documentclass[11pt,a4paper]{amsart}
\usepackage{amsmath,mathtools,multirow}

\usepackage[scr=boondox]{mathalpha}

\calclayout
\usepackage{hyperref}
\usepackage{anyfontsize}
\usepackage{stmaryrd}
\usepackage[shortlabels]{enumitem}
\usepackage{comment}
\usepackage{graphicx}
\usepackage{amsfonts}
\usepackage{amssymb}
\usepackage{wrapfig,lipsum,booktabs}
\usepackage{placeins}
\usepackage[hang]{subfigure}
\usepackage[utf8]{inputenc}
\usepackage[english]{babel}
\usepackage{amsthm}
\usepackage{subfigure} 
\usepackage[utf8]{inputenc}
\usepackage[margin=1in]{geometry} 
\usepackage[final]{pdfpages}
\newcounter{count}\numberwithin{count}{section}
\newtheorem{theorem}[count]{Theorem}
\newtheorem*{theorem*}{Theorem}

\newtheorem{lemma}[count]{Lemma}
\theoremstyle{definition}

\let\mathscr\relax
\usepackage{newpxtext,newpxmath}
\usepackage{blindtext}
\usepackage{setspace}
\author{Yazan Alamoudi}
\address{Department of Mathematics, University of Florida, Gainesville, FL 32611, United States}
\email{yazanalamoudi@ufl.edu}
\author{Krishnaswami Alladi}
\address{Department of Mathematics, University of Florida, Gainesville, FL 32611, United States}
\email{alladik@ufl.edu}

\title[\resizebox{15.2cm}{!}{Asymptotic estimates for sifted sums involving the Möbius function and the number of prime factors}]{Asymptotic estimates for sifted sums involving the Möbius function and the number of prime factors}

\begin{document}

\begin{abstract}
This paper gives refined quantitative estimates, valid in the range $1.9\le y\le x$, for $M_{2,\omega}(x,y)=\sum_{\substack{p_{1}(n)> y\\ n\leq x} } (\omega(n)-1)\mu(n)$ and $M_{\omega}(x,y)=\sum_{\substack{p_{1}(n)> y\\ n\leq x} } \omega(n)\mu(n)$ with $p_1(n)$ denoting the smallest prime factor. These are used to give sharp estimates for \[A_{2,\omega}(x)=\sup_{f\in\mathcal{B}}|\sum_{2\leq n\leq x}(\omega(n)-1)\mu(n)f(p(n))|\quad \text{and} \quad A_{\omega}(x)=\sup_{f\in\mathcal{B}}|\sum_{2\leq n\leq x}\omega(n)\mu(n)f(p(n))|,\] where $\mathcal{B}$ is the set of all functions absolutely bounded by $1$.
\end{abstract}
\maketitle
\section{introduction}
In the nineteen seventies, K. Alladi \cite{KA77} demonstrated a surprising duality between prime factors. Specifically, if $p_1(n)$ denotes the smallest prime factor of $n$ and $P_1(n)$ denotes the largest, then
\[\sum_{\substack{d|n\\d>1}}\mu(d)f(p_1(n))=-f(P_1(n))\quad \text{and} \quad \sum_{\substack{d|n\\d>1}}\mu(d)f(P_1(n))=-f(p_1(n)).\]
This was used to obtain the noteworthy\footnote{It is worthwhile to recall and compare the well-known equivalence between the prime number theorem and $\sum_{n\ge 1}\frac{\mu(n)}{n}=0$.} sum 
\[\sum_{\substack{p_1(n)\equiv r\pmod{m}\\n>1}}\frac{\mu(n)}{n}=\frac{-1}{\phi(m)}.\]
The above sum can be viewed as a special case of the following general correspondence on bounded functions
\[\sum_{n\leq x}f(P_1(n))\sim cx\Leftrightarrow\sum_{n\geq 2}\frac{\mu(n)f(p_1(n))}{n}=-c.\]
A key ingredient in the above correspondence is that ${\sum_{n\le x}\mu(n)f(p_1(n))=o(x)}$ for any bounded function $f$. The relation
\[\sum_{n\le x}\mu(n)f(p_1(n))=-\sum_{p\le x}f(p)M(\frac{x}{p},p),\]
with $M(x,y)=\sum_{\substack{p_{1}(n)> y\\ n\leq x} } \mu(n)$, suggests that refined estimates of $M(x,y)$ provide a useful tool for bounding the sum $\sum_{n\le x}\mu(n)f(p_1(n))$. In fact, for this reason, Alladi showed that $\frac{x\rho'(\alpha)}{\log y}$ is an asymptotic estimate for $M(x,y)$ when $\alpha>1$ is fixed.

Alladi also stated higher-order dualities involving $p_i(n)$, the $i$-th smallest prime factor of $n$, and $P_i(n)$ the $i$-th largest; the most relevant for this manuscript is \goodbreak
\[\sum_{\substack{d|n\\d>1}}\mu(d)(\omega(d)-1)f(p_1(n))=f(P_2(n)).\]
Given that
\begin{align*}
\sum_{2\leq n\leq x}(\omega(n)-1)\mu(n)f(p(n))&=-\sum_{p\leq x}f(p)M_{2,\omega}(\frac{x}{p},p)-\sum_{p\leq x}f(p)M(\frac{x}{p},p),\tag{1.1}\\
\intertext{and} 
\sum_{2\leq n\leq x}\omega(n)\mu(n)f(p(n))&=-\sum_{p\leq x}f(p)M_{\omega}(\frac{x}{p},p)-\sum_{p\leq x}f(p)M(\frac{x}{p},p),\tag{1.2}
\end{align*}
with $M_{2,\omega}(x,y)=\sum_{\substack{p_{1}(n)> y\\ n\leq x} } (\omega(n)-1)\mu(n)$ and $M_{\omega}(x,y)=\sum_{\substack{p_{1}(n)> y\\ n\leq x} } \omega(n)\mu(n)$, one is likewise motivated to give refined estimates of $M_{2,\omega}(x,y)$ and $M_{\omega}(x,y)$. Alamoudi \cite[Thm 1.1]{YAANT} has already found such estimates valid for ranges of the form $1.9\le y\le Y_0\exp(\mathscr{p}\frac{\log x}{(\log\log (x+1))^{1+\epsilon}})$. However, full range estimates have not been investigated in prior work. As such, the first main theorem of this manuscript will do just that. 
\begin{theorem} Let $m_\omega(\alpha)=H(\alpha-1)(-\int_{1}^{\alpha}\frac{\rho'(\alpha - t)}{t}dt-\frac{1}{\alpha})$ and $m_{2,\omega}(\alpha)=m_\omega(\alpha)-\rho'(\alpha)$, $\rho$ denotes the Dickman function and $H$ denotes the Heaviside function. Then, there are absolute constants $C$, $C'$, $C_{2,\omega}$, $C'_{2,\omega}$, $C'_{\omega}$, and $C_{\omega}$ such that whenever $1.9\le y\le x$ we have:
\begin{align*}
M_{2,\omega}(x,y)&=\frac{xm_{2,\omega}(\alpha)}{\log y}+ r_{2,\omega}(x,y),\tag{1.3}\\
M_{\omega}(x,y)\,\,&=\frac{xm_{\omega}(\alpha)}{\log y}+\frac{y}{\log y}+r_{\omega}(x,y),\tag{1.4}\\
\intertext{and}
M(x,y)\;&=\frac{x\rho'(\alpha)}{\log y}+\frac{y}{\log y}+r(x,y)\tag{1.5},
\end{align*}
with
\[\max\{\frac{|r_{2,\omega}(x,y)|}{C_{2,\omega}},\frac{|r_{\omega}(x,y)|}{C_{\omega}},\frac{|r(x,y)|}{C}\}\le\frac{x}{(\log y)^2},\]
and
\[\max\{\frac{|r_{2,\omega}(x,y)|}{C'_{2,\omega}},\frac{|r_{\omega}(x,y)|}{C'_{\omega}},\frac{|r(x,y)|}{C'}\}\le\frac{x\log\log(x+1)}{\sqrt[3]{(\log x)^2}}.\]

\end{theorem}
To better understand how Theorem 1.1 can be used, we note that Alladi and Johnson \cite[Thm~13]{KAJJ} obtained the following correspondences
\[\sum_{n\leq x}f(P_1(n))\sim\sum_{n\leq x}f(P_2(n))\sim \kappa x\Leftrightarrow\sum_{n\geq 2}\frac{\mu(n)\omega(n)f(p_1(n))}{n}=0,\]
under the assumption \cite[Thm~12]{KAJJ} that \[{\sum_{n\le x}\mu(n)\omega(n)f(p_1(n))=o(x)}.\]
In fact, Alladi--Johnson cite this paper in the proof of \cite[Thm~12]{KAJJ}; The proof of which follows from the stronger quantitative estimate below.
\begin{theorem}
Let $\mathcal{B}$ denote the set of functions absolutely bounded by $1$. Then
\[A_{2,\omega}(x)=\sup_{f\in\mathcal{B}}|\sum_{2\leq n\leq x}(\omega(n)-1)\mu(n)f(p(n))|=(2-\log(1+\sqrt{e})-2\rho(1+\sqrt{e}))\frac{2x}{\log x}+O(\frac{x(\log\log x)^2}{\sqrt{(\log x)^3}})\]
and 
\[A_{\omega}(x)=\sup_{f\in\mathcal{B}}|\sum_{2\leq n\leq x}\omega(n)\mu(n)f(p(n))|=(3-\log(1+\sqrt[3]{e^2})-3\rho(1+\sqrt[3]{e^2}))\frac{2x}{\log x}+O(\frac{x(\log\log x)^2}{\sqrt{(\log x)^3}}).\]
\end{theorem}
The above theorem verifies the correspondence but is in fact much stronger. Note that, just like the first order duality, setting $f$ equal to the indicator functions of integers congruent to $r$ modulo $m$ gives
\[\sum_{\substack{p_1(n)\equiv r\pmod{m}\\n>1}}\frac{(\omega(n)-1)\mu(n)}{n}=\frac{1}{\phi(m)},\]
which was proven unconditionally by Alladi--Johnson \cite[Eq.~9.2]{KAJJ} and later generalized by Alladi--Sengupta \cite{KASS}.
As we will see, the sifted sums $M_\omega$ and $M_{2,\omega}$ exhibit intriguing properties. For one, if $\alpha=1+\sqrt{e}$ then $M_\omega(x,x^{\frac{1}{\alpha}})\ll \frac{x}{(\log x)^2}$ but $M_\omega(x,x^{\frac{1}{\alpha}})\sim \frac{x\alpha m_{\omega}(\alpha)}{\log x}\asymp\frac{x}{\log x}$ for any fixed ${\alpha\in(1,\infty)\setminus\{1+\sqrt{e}\}}$. This is unique to $M_{\omega}(x,y)$ and does not occur for $M(x,y)$ because $m_\omega(\alpha)$ has a zero at $\alpha=1+\sqrt{e}$ and $\rho'$ is zero free on $(1,\infty)$. This, along with a more thorough overview of the analytic features of $m_\omega(\alpha)$ and $m_{2,\omega}(\alpha)$, will be the subject of the next section. Section 3 uses these analytic properties, together with Buchstab iteration, to establish Theorem 1.1.  Section 4 concludes the manuscript by utilizing the results of the prior two sections to obtain the sharp bounds on $A_\omega$ and $A_{2,\omega}$ stated in Theorem 1.2.
\section{Analytic preliminaries}
The following theorem connects $m_\omega(\alpha)$ and $m_{2,\omega}(\alpha)$ to a continuous analog of Buchstab iteration.\footnote{On January 14, 2025, K. Alladi sent Y. Alamoudi the delay differential equation for $m_\omega(\alpha)$, though the notes themselves were dated April 17, 2023. K. Alladi was able to find it by analogy with Buchstab iteration adapted to a continuous setting. Alamoudi found the solution $m_\omega(\alpha)=H(\alpha-1)(-\int_{1}^{\alpha}\frac{\rho'(\alpha - t)}{t}dt-\frac{1}{\alpha})$ via Laplace transforms and sent it to Alladi on January 18, 2025. Shortly after, Alamoudi presented this in his oral exam \cite{YAORAl}.}
\begin{theorem}
Let $m_{2,\omega}(\alpha)$ and $m_{\omega}(\alpha)$ be defined by $m_{2,\omega}(\alpha)=0$
and $m_{\omega}(\alpha)=m(\alpha)=\rho'(\alpha)$ when $0<\alpha \leq2$, and
\[
   \{\alpha m_{2,\omega}(\alpha)\}'=-m_{2,\omega}(\alpha-1)-m(\alpha-1)\quad \text{and} \quad \{\alpha m_{\omega}(\alpha)\}'=-m_{\omega}(\alpha-1)-m(\alpha-1)
\]
when $\alpha>2$.
Then $m_{2,\omega}(\alpha)=m_\omega(\alpha)-\rho'(\alpha)$ and $m_\omega(\alpha)=H(\alpha-1)(-\int_{1}^{\alpha}\frac{\rho'(\alpha - t)}{t}dt-\frac{1}{\alpha})$.
\end{theorem}

\begin{proof}
Subtracting $\{\alpha m(\alpha)\}'=m(\alpha-1)$, where $m$ denotes $\rho'$, from the delay differential equation for $m_\omega(\alpha)$ gives a solution that vanishes for $\alpha\leq2$ and so $m_{2,\omega}(\alpha)=m_\omega(\alpha)-\rho'(\alpha)$. Thus, it suffices to solve the delay differential equation for $m_{2,\omega}$. Taking Laplace transforms of the equation for $m_{2,\omega}$ gives\goodbreak 
\[\frac{d}{ds}\mathcal{L}(m_{2,\omega})-\frac{e^{-s}}{s}\mathcal{L}(m_{2,\omega})=\frac{e^{-s}}{s}\mathcal{L}(m).\]
Using the integrating factor, $I(s)=e^{\int_s^\infty\frac{e^{-t}}{t}}$ we can write
\[\mathcal{L}(m_{2,\omega})=\frac{-1}{I(s)}\left(\int_{s}^\infty I(s)\frac{e^{-s}}{s}\mathcal{L}(m)ds+C_0\right).\]
Noting that the well-known formula for the Laplace transform of the Dickman function \cite[Thm 7.10]{MV} can be rewritten as \(\mathcal{L}(\rho)=\frac{e^{-\int_s^\infty\frac{e^{-t}}{t}}}{s}=\frac{1}{sI(s)}\), and simplifying the above becomes
\[\mathcal{L}(m_{2,\omega})=\mathcal{L}(m)\int_{s}^\infty-\frac{e^{-t}}{t}dt+C_0\mathcal{L}(m)+\int_{s}^\infty-\frac{e^{-t}}{t}dt+C_0+1.\]
Taking inverse Laplace transforms in the above gives
\[m_{2,\omega}(\alpha)=H(\alpha-1)(-\int_{1}^{\alpha}\frac{\rho'(\alpha - t)}{t}dt-\frac{1}{\alpha})+C_0\rho'(\alpha) +(C_0+1)\delta(\alpha).\]
The initial conditions force $C_0=-1$, and the conclusion follows. 
\end{proof}
The following lemma establishes properties of $m_\omega(\alpha)$ and $m_{2,\omega}(\alpha)$ that will be of use in the subsequent sections.
\begin{lemma}
The following assertions hold true.
\begin{enumerate}[label=(\roman*)]
\item For $\alpha>0$, $m_\omega(\alpha)+\frac{10}{\alpha}$ is monotone decreasing,\\
\item $m_\omega(\alpha)=\frac{e^\gamma}{\alpha^2}+O(\frac{1}{\alpha^3})$,\\
\item $m_{2,\omega}(\alpha)=\frac{e^\gamma}{\alpha^2}+O(\frac{1}{\alpha^3})$,\\
\item and for all $\alpha\geq3$ we have $ m_\omega(\alpha)+\rho'(\alpha)\geq0$.\\
\end{enumerate}
\end{lemma}

\begin{proof}
For (i), the case where $\alpha \le2$ follows from the initial condition $m_\omega(\alpha)=\rho'(\alpha)$. When $\alpha>2$ the proof follows from
\[m_\omega(\alpha)+\frac{10}{\alpha}=\int_{1}^{\alpha-1}\frac{\rho(\alpha - t)}{t^2}dt+\frac{1}{\alpha-1}+\frac{9}{\alpha}-\rho(\alpha-1)\]
\[\Rightarrow \frac{d}{d\alpha}(m_\omega(\alpha)+\frac{10}{\alpha})=\int_{1}^{\alpha-1}\frac{\rho'(\alpha - t)}{t^2}dt-\frac{1}{(\alpha-1)^2}-\frac{9}{\alpha^2}-\rho'(\alpha-1)<0.\]
To prove (ii) and (iii), begin with
\[P(t)=\int_0^t\rho(t) dt=e^\gamma-\int_t^\infty\rho(t) dt,\]
where the second equality follows from the well-known $\int_0^\infty\rho(t) dt=e^\gamma$. Now 
\[\int_t^\infty\rho(t) dt\le\int_t^\infty\frac{1}{\Gamma(t+1)} dt\ll\frac{1}{\Gamma(t-1)}\int_t^\infty\frac{1}{t(t-1)} dt\ll\frac{1}{\Gamma(t)}\]
hence
\[P(t)=e^\gamma+O(\frac{1}{\Gamma (t)}).\]
\goodbreak
On the other hand, by integration by parts\smallskip
\begin{align*}
m_\omega(\alpha)&=\int_{1}^{\alpha -1}\frac{\rho(\alpha -t)}{t^2}dt+\frac{1}{\alpha(\alpha-1)}+O(\frac{1}{\alpha^3})\\
&=P(\alpha-1)+2\int_{1}^{\alpha -1}\frac{P(\alpha -t)}{t^3}dt+O(\frac{1}{\alpha^3})\\
&=e^\gamma+2\int_{1}^{\alpha-1}\frac{e^\gamma}{t^3}dt+O(\int_{1}^{\alpha -1}\frac{1}{\Gamma(\alpha -t)t^3}dt)+O(\frac{1}{\alpha^3})\\
&=\frac{e^\gamma}{\alpha^2}+O(\int_{1}^{\alpha -1}\frac{1}{\Gamma(\alpha -t)t^3}dt)+O(\frac{1}{\alpha^3}).
\end{align*}
The conclusion for (ii) follows once we note that
\[\int_{1}^{\alpha -1}\frac{1}{\Gamma(\alpha -t)t^3}dt\ll\int_{1}^{\frac{\alpha}{2}}\frac{1}{\Gamma(\alpha -t)t^3}dt+\int_{\frac{\alpha}{2}}^{\alpha}\frac{1}{\Gamma(\alpha -t)t^3}dt\ll\frac{1}{\Gamma(\frac{\alpha}{2})}+\frac{1}{\alpha^3}\ll\frac{1}{\alpha^3},\]
and (iii) follows from $m_{2,\omega}(\alpha)=m_{\omega}(\alpha)-\rho'(\alpha)$.

Lastly, to prove (iv), start with
\[ m_\omega(\alpha)+\rho'(\alpha)=P(\alpha-1)+2\int_{1}^{\alpha -1}\frac{P(\alpha -t)}{t^3}dt-\frac{1}{(\alpha-1)^2\alpha}-\rho(\alpha-1)+\rho'(\alpha).\]
The proof immediately follows since for $\alpha\ge3$ the above is
\[\geq 2+e^\gamma-\frac{2\alpha+1}{\alpha}\left(\frac{1}{(\alpha-1)^2}+\frac{1}{\Gamma(\alpha)}\right)\ge2+e^\gamma-\frac{7}{4}\ge0.\]
\end{proof}
The analytic properties listed so far will be sufficient for the goals of this manuscript. In particular, the next section will utilize the delay differential equation and the first three parts of Lemma 2.2 to establish Theorem 1.1.
\section{Proof of Theorem 1.1}
The proof of Theorem 1.1 will utilize Buchstab iteration and Stieltjes integration. Observe that choosing $f$ to be the indicator function for primes greater than $\sqrt{x}$ in (1.1) and (1.2) gives
\enlargethispage{0.5\baselineskip}
\begin{align*}
 M_{\omega}(x,y)\;&=M_{\omega}(x,\sqrt x)-\sum_{y<p\leq \sqrt x } \{M_\omega(\frac{x}{p},p)+M(\frac{x}{p},p)\},\tag{3.1}\\
M_{2,\omega}(x,y)&=M_{2,\omega}(x,\sqrt x)-\sum_{y<p\leq \sqrt x } \{M_{2,\omega}(\frac{x}{p},p)+M(\frac{x}{p},p)\};\tag{3.2}\\ 
\intertext{moreover,}
M(x,y)\;&=M(x,\sqrt x)-\sum_{y<p\leq \sqrt x } M(\frac{x}{p},p)\tag{3.3}
\end{align*}
can similarly be obtained from $\sum_{n\le x}\mu(n)f(p_1(n))=-\sum_{p\le x}f(p)M(\frac{x}{p},p)$ and was already noted in more general form in \cite{KA82}. On the other hand, the discrepancy from Stieltjes integration will be controlled by the lemma below.

\begin{lemma}
Suppose $f_x(t)$ is differentiable absolutely bounded monotone decreasing function of $t$, and $B\in C^1$ with $|R(t)|\leq B(t)$, then
\[|\int_{a}^{b}f_x(t)dR(t)|\le3\max(|f_x(a)|,|f_x(b)|)\sup_{[a,b]}|B(t)|.\tag{3.4}\]
\end{lemma}
\begin{proof}
Since $f'_x$ does not change signs, an application of Stieltjes integration leads to
\[|\int_{a}^{b}f_x(t)dR(t)|\le|f_x(t)R(t)|_{a}^{b}|+|f_x(t)B(t)|_{a}^{b}|+|\int_{a}^{b}f_x(t)B'(t)dt|.\]
At this point, note that the given conditions on $f_x$ allow the use of the second mean value theorem of integrals (possibly after a shift by a constant) to get that the above is 
\[\le|f_x(t)R(t)|_{a}^{b}|+|f_x(t)B(t)|_{a}^{b}|+\sup_{[a,b]}|f_x|\sup_{\eta\in[a,b]}|\int_{a}^{\eta}B'(t)dt|\le3\max(|f_x(a)|,|f_x(b)|)\sup_{[a,b]}|B(t)|.\]

\end{proof}
We will first refine the bound for $M(x,y)$ (see \cite{KA82}) as it will be used in the induction proof of $M_{2,\omega}(x,y)$.

\begin{proof}[Proof of (1.5) with $y$-dependent error]
Induction on $\lceil\alpha\rceil$. Base case $\alpha \le 2$ follows from
$$|M(x,y)+\pi(x)|\leq|M(x,\sqrt{x})+\pi(x)|\leq \pi(\sqrt x)+1$$
 $$\implies\max\{|M(x,y)+\frac{x}{\log x}|,|M(x,\sqrt x)+\frac{x}{\log x}|\}\leq \pi(\sqrt x)+1+|\pi(x)-\frac{x}{\log x}|\ll \frac{x}{(\log x)^2}.$$
Suppose this is true up to $\lceil\alpha\rceil=N$. Then, summing over $p$ and Stieltjes integration give
\[\hspace*{-0.7cm}|\sum_{y<p\le \sqrt x}M(\frac{x}{p},p)-\int_y^{x^\frac{1}{2}}\frac{x\rho'(\frac{\log x}{\log t}-1)}{t(\log t)^2}dt|\leq|\int_y^{x^\frac{1}{2}}\frac{x\rho'(\frac{\log x}{\log t}-1)}{t\log t}dR(t)|+ \sum_{y<p\le \sqrt x}\frac{Cx}{p(\log p)^2}+\sum_{p\le \sqrt x}\frac{p}{\log p}.\]
On the other hand,\[|\sum_{y<p\le \sqrt x}M(\frac{x}{p},p)-\int_y^{x^\frac{1}{2}}\frac{x\rho'(\frac{\log x}{\log t}-1)}{t(\log t)^2}dt|=|\frac{x}{\log x}+\sum_{y<p\le \sqrt x}M(\frac{x}{p},p)+\frac{x\rho'(\alpha)}{\log y}|,\]
which when combined with $|M(x,\sqrt x)+\frac{x}{\log x}|\ll \frac{x}{(\log x)^2}$ and (3.3) gives
\[\mathclap{|M(x,y)-\frac{x\rho'(\alpha)}{\log y}|\leq|\int_y^{x^\frac{1}{2}}\frac{x\rho'(\frac{\log x}{\log t}-1)}{t\log t}dR(t)|+|\int_y^{x^\frac{1}{2}}\frac{Cx}{t(\log t)^2}dR(t)|+\int_y^{x^\frac{1}{2}}\frac{Cx}{t(\log t)^3}dt+C''\frac{x}{(\log x)^2}.}\]
Applying Lemma 3.1 on the first two integrals \footnote{To meet the conditions of Lemma 3.1 for the first integral, write $\int_y^{x^\frac{1}{2}}\frac{x\rho'(\frac{\log x}{\log t}-1)}{t\log t}dR(t)|\leq |\int_y^{x^\frac{1}{2}}\frac{x(1-\rho'(\frac{\log x}{\log t}-1))}{t\log t}dR(t)|+|\int_y^{x^\frac{1}{2}}\frac{x}{t\log t}dR(t)|$ and apply to each.} with $B(t)=\frac{A}{\exp{\sqrt{c\log t}}}$  and simplifying gives
\[\leq\frac{x}{(\log y)^2}(\frac{6A\log y}{\exp{ \sqrt{c\log{ y}}}}+\frac{A3C}{\exp{ \sqrt{c\log{ y}}}}+\frac{C}{2})+(C''-2C)\frac{x}{(\log x)^2}.\]
The conclusion follows by picking $y>Y_0$ sufficiently large and $C>C''$. The finitely many cases where $y\leq Y_0$ are handled by Alamoudi's result via the Selberg--Delange method. Alternatively, one could replace $\exp{ \sqrt{c\log{ y}}}$ in $B(y)$ with $\exp{ \sqrt{ L\log{ y}}}$ for some sufficiently large $L$, which is allowed in view of the Chudakov or Vinogradov zero-free regions (see\cite{MV,Ten}), so that
\[(\frac{6A\log y}{\exp{ \sqrt{L \log{ y}}}}+\frac{A3C}{\exp{ \sqrt{ L{ \log{ y}}}}}+\frac{C}{2})<C\]
for $y\geq2$. In any case, the first claim is proven. 
\end{proof}
We can now complete the proof of the quantitative statement Theorem 1.1 with the $y$-dependent error bound. This is done by repeating the above process for $M_{2,\omega}$, which will immediately give the result for $M_\omega$ by taking the difference.
\begin{proof}[Proof of (1.3) and (1.4) with $y$-dependent error]
This is done by induction on $\lceil\alpha\rceil$. The base case is $\alpha \le 2$ is $M_{2,\omega}(x,y)=M_{2,\omega}(x,\sqrt{x})=1$. 
Suppose this is true up to $\lceil\alpha\rceil=N$. Then,
\[|\sum_{y<p\le \sqrt x}M_{2,\omega}(\frac{x}{p},p)-\sum_{y<p\le \sqrt x}\frac{xm_{2,\omega}(\frac{\log x}{\log p}-1)}{p\log p}|\leq \sum_{y<p\le \sqrt x}\frac{C_{2,\omega}x}{p(\log p)^2}.\]
Because the Buchstab relation for $M_{2,\omega}(x,y)$ involves $M(x,y)$, the corresponding sum must be added; this, combined with Stieljets integration gives
\[|\sum_{y<p\le \sqrt x}\{M_{2,\omega}(\frac{x}{p},p)-M(\frac{x}{p},p)\}-\int_y^{x^\frac{1}{2}}\frac{x\{m_{2,\omega}(\frac{\log x}{\log t}-1)-\rho'(\frac{\log x}{\log t}-1)\}}{t(\log t)^2}dt|\]
\[\leq|\int_y^{x^\frac{1}{2}}\frac{x\rho'(\frac{\log x}{\log t}-1)}{t\log t}dR(t)|+|\int_y^{x^\frac{1}{2}}\frac{xm_{2,\omega}(\frac{\log x}{\log t}-1)}{t\log t}dR(t)|+ \sum_{y<p\le \sqrt x}\frac{(C+C_{2,\omega})x}{p(\log p)^2}+\sum_{p\le \sqrt x}\frac{p}{\log p}.\]
Now, to meet the conditions of Lemma 3.1, use Lemma 2.2, write \[\int_y^{x^\frac{1}{2}}\frac{xm_\omega(\frac{\log x}{\log t}-1)}{t\log t}dR(t)|\leq |\int_y^{x^\frac{1}{2}}\frac{x(\frac{10}{(\frac{\log t}{\log x}-1)}+m_\omega(\frac{\log x}{\log t}-1))}{t\log t}dR(t)|+|\int_y^{x^\frac{1}{2}}\frac{x\frac{10}{(\frac{\log t}{\log x}-1)}}{t\log t}dR(t)|,\] and apply to each. This gives,
\[\int_y^{x^\frac{1}{2}}\frac{xm_\omega(\frac{\log x}{\log t}-1)}{t\log t}dR(t)|\leq 3A\frac{x\mathscr{M}}{\log y\exp{ \sqrt{c\log{ y}}}},\]
where $\sup_{\alpha>0}\{|m_\omega(\alpha)|\}+20\leq\mathscr{M}<\infty$.
Combining the above with the computation from the previous case and simplifying gives\vspace*{0.35cm}
\[\vspace*{0.35cm}\mathclap{|M_{2,\omega}(x,y)-\frac{xm_{2,\omega}(\alpha)}{\log y}|\leq\frac{x}{(\log y)^2}(\frac{3A(\mathscr{M}+4)\log y}{\exp{ \sqrt{c\log{ y}}}}+\frac{A3(C+C_{2,\omega})}{\exp{ \sqrt{c\log{ y}}}}+\frac{(C+C_{2,\omega})}{2})+(C'-2(C+C_{2,\omega}))\frac{x}{(\log x)^2}+1.}\]

The conclusion follows by the same line of reasoning as in the $\rho'$ case (and noting that the ''$+1$'' at the end could be trivially accounted for by adjusting the constants appropriately).
\end{proof}
The proof of Theorem 1.1 is complete once we establish the $y$-independent error bounds.
\begin{proof}[Proof of y-independent error bound]
By \cite[Thm 1.1]{YAANT} for $1.9\le y\le\exp((\log x)^{\frac{3}{4}})$ we have absolute constants $C_{1,1},C_{1,0}$ such that \goodbreak
\[M_{2,\omega}(x,y)=\frac{e^\gamma x}{\alpha^2 \log y}+\varepsilon_{1,1}(x,y)\quad \text{and}\quad M(x,y)=\varepsilon_{1,0}(x,y),\]
with \[\max{(\frac{|\varepsilon_{1,1}(x,y)|}{C_{1,1}},\frac{|\varepsilon_{1,0}(x,y)|\log \log(x+1)}{C_{1,0}})}\leq \frac{x\log\log (x+1)}{\log x}\bigg\{(\frac{\log y}{\log x})^{2}+\frac{1}{\log x\exp (c''\sqrt{\log y})}\bigg\}.\] This directly gives $|\varepsilon_{1,0}(x,y)|\leq C''\frac{x\log\log(x+1)}{\sqrt{(\log x)^3}}$ for some absolute constant $C''$. To get $M_{2,\omega}$ in the form of Theorem 1.1, we utilize the lemma from the previous section to get
\[M_{2,\omega}(x,y)=\frac{xm_{2,\omega}(\alpha)}{\log y}+(\frac{e^\gamma}{\alpha^2}-m_{2,\omega}(\alpha))\frac{x}{\log y}+\varepsilon_{1,1}(x,y)=\frac{xm_{2,\omega}(\alpha)}{\log y}+ r'_{2,\omega}(x,y)\]
with 
\[|r'_{2,\omega}(x,y)|\leq|\frac{e^\gamma}{\alpha^2}-m_{2,\omega}(\alpha)|\frac{x}{\log y}+|\varepsilon_{1,1}(x,y)|\leq C''_{2,\omega}\frac{x\log\log(x+1)}{\sqrt{(\log x)^3}}.\] On the other hand, $y\ge\exp((\log x)^{\frac{3}{4}})$ the bound for both $M$ and $M_{2,\omega}$ is a consequence of \({\frac{x}{(\log y)^2}\leq \frac{x}{(\log x)^\frac{3}{2}}}\) and the $y$ dependent error. The conclusion follows as $C'=\max(C'',C)$ and ${C'_{2,\omega}=\max(C''_{2,\omega},C_{2,\omega})}$ satisfy the claim of the theorem.
\end{proof}
\section{Proof of theorem 1.2}
All that is left is Theorem 2.1. To start, it is convenient to note that (1.1) could be simplified and written as
\[\sum_{2\leq n\leq x}(\omega(n)-1)\mu(n)f(p(n))=-\sum_{p\leq x}f(p)M_{\omega}(\frac{x}{p},p).\tag{4.1}\]
From (4.1) and (1.2) it can be seen that
\[A_{2,\omega}(x)=\sum_{p\leq x}|M_{\omega}(\frac{x}{p},p)|\quad \text{and}  \quad A_{\omega}(x)=\sum_{p\leq x}|M_{\omega}(\frac{x}{p},p)+M(\frac{x}{p},p)|.\]
To begin with, note that $\sum_{\sqrt x\leq p\leq x}|M_{\omega}(\frac{x}{p},p)|<<1$ and $\sum_{\sqrt x\leq p\leq x}|M(\frac{x}{p},p)|=\frac{x}{\log x}+O(\frac{x}{(\log x)^2})$ and so
\[A_{2,\omega}(x)=\sum_{p\leq \sqrt{x}}|M_{\omega}(\frac{x}{p},p)|+O(1) \quad \text{and} \quad A_{\omega}(x)=\sum_{p\leq \sqrt{x}}|M_{\omega}(\frac{x}{p},p)+M(\frac{x}{p},p)|+\frac{x}{\log x}+O(\frac{x}{(\log x)^2}).\]
We complete the argument by Stieltjes integration on the quantitative statement with a $ y$-independent error term starting with $A_{2,\omega}(x)$. This gives, \smallskip
\begin{align*}
A_{2,\omega}(x)&=\int_{2}^{\sqrt{x}}\frac{x|m_{\omega}(\frac{\log x}{\log t}-1)|}{t(\log t)^2}dt+O(\int_{2}^{\sqrt{x}}\frac{x\log\log(\frac{x}{t}+1)}{t\sqrt{(\log \frac{x}{t})^3}\log t}dt)\\
&=\int_{2}^{\sqrt{x}}\frac{x|m_{\omega}(\frac{\log x}{\log t}-1)|}{t(\log t)^2}dt+O(\frac{x(\log\log x)^2}{\sqrt{(\log x)^3}})\\
&=\frac{x}{\log x}\int_{2}^{\frac{\log x}{\log 2}}|m_{\omega}(u-1)|du+O(\frac{x(\log\log x)^2}{\sqrt{(\log x)^3}})\tag{4.2}.
\end{align*}
Now, although Lemma 2.2 shows that $m_\omega(\alpha)\geq m_w(\alpha)+\rho'(\alpha)\geq0$ when $\alpha\geq 3$, both are initially negative. This can be seen from the explicit formulas 
\[m_\omega(\alpha)=-\frac{1}{\alpha}+\frac{2}{\alpha}\log(\alpha-1)=\frac{1}{\alpha}+2\rho'(\alpha)\tag{4.3}\]
and 
\[m_\omega(\alpha)+\rho'(\alpha)=-\frac{2}{\alpha}+\frac{3}{\alpha}\log(\alpha-1)=\frac{1}{\alpha}+3\rho'(\alpha),\tag{4.4}\]
valid when $2\le\alpha\le3$. In particular, $m_\omega(\alpha)$ flips signs at $1+\sqrt{e}$ and $m_\omega(\alpha)+\rho'(\alpha)$ does so at $1+\sqrt[3]{e^2}$.
Hence, (4.2) becomes
\medskip
\begin{align*}
A_{2,\omega}(x)&=\frac{x}{\log x}\left(\int_{2}^{\frac{\log x}{\log 2}}m_{\omega}(u-1)du-2\int_{2}^{2+\sqrt{e}}m_{\omega}(u-1)du\right)+O(\frac{x(\log\log x)^2}{\sqrt{(\log x)^3}})\\
&=\frac{x}{\log x}\left(\int_{2}^{\infty}m_{\omega}(u-1)du-2\int_{2}^{2+\sqrt{e}}(\frac{1}{u-1}+2\rho'(u-1))du\right)+O(\frac{x(\log\log x)^2}{\sqrt{(\log x)^3}})\\
&=\frac{x}{\log x}\left(-2\int_{2}^{2+\sqrt{e}}(\frac{1}{u-1}+2\rho'(u-1))du\right)+O(\frac{x(\log\log x)^2}{\sqrt{(\log x)^3}})\\
&=(2-\log(1+\sqrt{e})-2\rho(1+\sqrt{e}))\frac{2x}{\log x}+O(\frac{x(\log\log x)^2}{\sqrt{(\log x)^3}}).
\end{align*}
\medskip

In exactly the same way, by Stieltjes integration on the quantitative statement with the $y$-independent error term for $A_{\omega}(x)$ we have
\medskip
\begin{align*}
\hspace{-0.5cm}
A_{\omega}(x)&=\frac{x}{\log x}(1+\int_{2}^{\frac{\log x}{\log 2}}|m_{\omega}(u-1)+\rho'(u-1)|du)+O(\frac{x(\log\log x)^2}{\sqrt{(\log x)^3}})\\
&=\frac{x}{\log x}\left(1+\int_{2}^{\frac{\log x}{\log 2}}m_{\omega}(u-1)+\rho'(u-1)du-2\int_{2}^{2+\sqrt[3]{e^2}}m_{\omega}(u-1)+\rho'(u-1)du\right)+O(\frac{x(\log\log x)^2}{\sqrt{(\log x)^3}})\\
&=\frac{x}{\log x}\left(1+\int_{2}^{\infty}m_{\omega}(u-1)+\rho'(u-1)du-2\int_{2}^{2+\sqrt[3]{e^2}}\frac{1}{u-1}+3\rho'(u-1)du\right)+O(\frac{x(\log\log x)^2}{\sqrt{(\log x)^3}})\\
&=\frac{x}{\log x}\left(-2\int_{2}^{2+\sqrt[3]{e^2}}\frac{1}{u-1}+3\rho'(u-1)du\right)+O(\frac{x(\log\log x)^2}{\sqrt{(\log x)^3}})\\
&=(3-\log(1+\sqrt[3]{e^2})-3\rho(1+\sqrt[3]{e^2}))\frac{2x}{\log x}+O(\frac{x(\log\log x)^2}{\sqrt{(\log x)^3}}).
\end{align*}
\medskip
This completes the proof of Theorem 1.2.

\bibliographystyle{amsplain}

\end{document}